\documentstyle[11pt,fleqn]{article}
\def\beq{\begin{equation}}
\def\eeq{\end{equation}}
\def\beqn{\begin{eqnarray}}
\def\eeqn{\end{eqnarray}}
\def\ni{\noindent}

\def\Uqso{U'_q({\rm so}_n)}

\begin{document}

%\begin{center}
\ni {\LARGE \bf
On the Casimir Elements of $q$-Algebras $\Uqso$\\
and Their Eigenvalues in Representations}
%\end{center}

\vspace{6mm}
\ni
%\bigskip
%\bigskip
{\sl A. M.~GAVRILIK and N. Z.~IORGOV}\\

\vspace{.2mm}
\ni
{\it Bogolyubov Institute for Theoretical Physics,  \\
Metrologichna Street 14b,  Kyiv 143, \ Ukraine}
%\bigskip

\vspace{1mm}
\begin{abstract}
\small
The nonstandard $q$-deformed algebras $\Uqso$ are known to possess
$q$-analogues of Gel'fand--Tsetlin type representations. For these
$q$-algebras, all the Casimir elements (corresponding
to basis set of Casimir elements of ${\rm so}_n$) are found,
and their eigenvalues within irreducible representations are given explicitly.
\end{abstract}

\vspace{1mm}
%\medskip
\ni{\large\bf 1. Introduction.}

\vspace{.2mm}
%\smallskip
The nonstandard deformation $\Uqso$, see \cite{GK}, of the Lie algebra
${\rm so}_n$ admits, in contrast to standard deformation \cite{DJ}
of Drinfeld and Jimbo, an explicit construction
of irreducible representations \cite{GK,GI} corresponding to those of
Lie algebra ${\rm so}_n$ in Gel'fand--Tsetlin formalism.
Besides, as it was shown in
\cite{NUW}, $\Uqso$ is the proper dual for the standard $q$-algebra
$U_q({\rm sl}_2)$ in the $q$-analogue of dual pair $({\rm so}_n,{\rm sl}_2)$.

Let us mention that
the algebra $U'_q({\rm so}_3)$ appeared earlier in the papers \cite{FO}.
As a matter of interest, this algebra arose naturally as the algebra
of observables \cite{NRZ} in 2+1 quantum gravity with 2D space fixed as torus.
At $n>3$, the algebras $\Uqso$ are no less important, serving
as intermediate algebras in deriving the algebra of observables in
2+1 quantum gravity with 2D space of genus $g>1$, so that $n$ depends on $g$,
$n=2g+2$ \cite{NR,G}. In order to obtain the algebra of observables,
the $q$-deformed algebra $U'_q({\rm so}_{2g+2})$ should be quotiented by
some ideal generated by (combinations of) Casimir elements of this algebra.
This fact, along with others, motivates the study of Casimir elements of
$\Uqso$.

%\newpage
\vspace{.4mm}
\ni {\large\bf 2. The $q$-deformed agebras $\Uqso$.}

\vspace{.2mm}
According to \cite{GK}, the nonstandard $q$-deformation $\Uqso$
of the Lie  algebra ${\rm so}_n$ is
given as a complex associative algebra with ${n-1}$
generating elements $I_{21}$, $I_{32}, \ldots$, $I_{n,n-1}$
obeying the defining relations (denote $q+q^{-1}\equiv [2]_q$)
\beq
\begin{array}{l}
I_{j,j-1}^2I_{j-1,j-2} + I_{j-1,j-2}I_{j,j-1}^2 -
[2]_q \ I_{j,j-1}I_{j-1,j-2}I_{j,j-1} = -I_{j-1,j-2}, \\[2mm]
I_{j-1,j-2}^2I_{j,j-1} + I_{j,j-1}I_{j-1,j-2}^2 -
[2]_q \ I_{j-1,j-2}I_{j,j-1}I_{j-1,j-2} = -I_{j,j-1},  \\[2mm]
[I_{i,i-1},I_{j,j-1}] =0  \qquad {\rm if} \quad \mid {i-j}\mid >1.
\end{array}  \label{f1}
\eeq
%%%%%%%%%%%%%%%%%%%%%%%%%%%%%%%%%%%%%%%%%%%%%%%%
Along with definition in terms of trilinear relations,
we also give a {\it `bilinear' presentation}. To this end, one
introduces the generators (here $k > l+1, \ \  1\leq k,l \leq n$)
\[
I^{\pm}_{k,l}\equiv [I_{l+1,l} , I^{\pm}_{k,l+1}]_{q^{\pm 1}}
\equiv q^{\pm 1/2}I_{l+1,l} I^\pm_{k,l+1} -
q^{\mp 1/2}I^\pm_{k,l+1} I_{l+1,l}
\]
together with $I_{k+1,k} \equiv I^+_{k+1,k}\equiv I^-_{k+1,k}$.
Then (\ref{f1}) imply
\[
[I^+_{lm} , I^+_{kl}]_q = I^+_{km}, \ \ \
[I^+_{kl} , I^+_{km}]_q = I^+_{lm}, \ \ \
[I^+_{km} , I^+_{lm}]_q = I^+_{kl}  \ \ \ {\rm if} \ \ k>l>m,
\]
\beq                                                \label{f4}
[I^+_{kl} , I^+_{mp}] = 0 \ \ \ \ \ \ \ \ {\rm if} \ \ \ \ k>l>m>p
\ \ \ \ \ \ {\rm or} \ \ \ \  k>m>p>l;
\eeq
%\vspace{2mm}
\[
[I^+_{kl} , I^+_{mp}] = (q-q^{-1}) (I^+_{lp}I^+_{km}-
I^+_{kp}I^+_{ml}) \ \ \ \ \ \ {\rm if} \ \ k>m>l>p.
\]
%%%%%%%%
Analogous set of relations exists involving $I_{kl}^-$ along with
$q\to q^{-1}$ (let us denote this ``dual'' set by (\ref{f4}$'$)). If $q\to 1$
(`classical' limit), both (\ref{f4}) and (\ref{f4}$'$)
reduce to those of ${\rm so}_n$.

Let us give explicitly the two examples, namely $n=3$
(the Odesskii--Fairlie algebra \cite{FO}) and $n=4$, using
the definition $[X,Y]_q\equiv q^{1/2}X Y - q^{-1/2} Y X$:
\[ \hspace{-0.2cm}
U'_q({\rm so}_4): \hspace{0.2cm}
\left\{ \hspace{0.2cm}
\parbox{12cm}
{
\beq
\hspace{-.8cm}
U'_q({\rm so}_3): \hspace{0.6cm}
[I_{21} , I_{32}]_q = I_{31}^+, \ \
[I_{32} , I_{31}^+]_q = I_{21}, \ \
[I_{31}^+ , I_{21}]_q = I_{32}.                       \label{FO}
\eeq
\beq
\hspace{-1.0cm}
\begin{array}{cclll}

 \hspace{0.4cm}& \hspace{0cm}&
[I_{32},I_{43}]_q = I_{42}^+,\ \ &
[I_{31}^+,I_{43}]_q = I_{41}^+,\ \ &
[I_{21},I_{42}^+]_q = I_{41}^+, \\
&&
[I_{43},I_{42}^+]_q=I_{32},&
[I_{43},I_{41}^+]_q=I_{31}^+,&
[I_{42}^+,I_{41}^+]_q=I_{21}, \\
&&
[I_{42}^+,I_{32}]_q=I_{43},&
[I_{41}^+,I_{31}^+]_q=I_{43},&
[I_{41}^+,I_{21}]_q=I_{42}^+,
\end{array}                                                      \label{O4}
\eeq
\vspace{-0.4cm}

\beq
\hspace{-1.0cm}
[I_{43},I_{21}]=0,\ \ \ [I_{32},I_{41}^+]=0,\ \ \
[I_{42}^+,I_{31}^+]=(q-q^{-1})(I_{21}I_{43}-I_{32}I_{41}^+).     \label{O4p}
\eeq
}
\right.
\]

The first relation in (\ref{FO}) can be viewed as the definition for
third generator
needed to give the algebra in terms of $q$-commutators.
Dual copy of the algebra $U'_q({\rm so}_3)$ involves the generator
$I_{31}^-=[I_{21},I_{32}]_{q^{-1}}$ and the other
two relations similar to (\ref{FO}), but with $q\to q^{-1}$.

In order to describe the basis of $\Uqso$  we introduce a
lexicographical ordering for the elements $I_{k,l}^+$ of $\Uqso$ with
respect to their indices, i.e., we suppose that $I_{k,l}^+\prec I_{m,n}^+$
if either $k<m$, or both $k=m$ and $l<n$. We define {\it an ordered
monomial} as the product of non--decreasing sequence of elements
$I_{k,l}^+$ with different $k,l$ such that $1\le l<k\le n$. The
following proposition describes the
Poincar${\acute {\rm e}}$--Birkhoff--Witt basis for the
algebra $\Uqso$.

\ni {\bf Proposition.} {\it The set of all ordered monomials is a basis
of $\Uqso$.}

\medskip
\ni {\bf 3. The Casimir elements of $\Uqso$.}
\medskip

As it is well-known, tensor operators of Lie algebras ${\rm so}_n$
are very useful in construction of invariants of these algebras.
With this in mind, let us introduce $q$-analogues of tensor operators
for the algebras $\Uqso$ as follows:
\beq
J^\pm_{k_1,k_2,\ldots,k_{2r}}=q^{\mp\frac{r(r-1)}{2}}
\mathop{{\sum}'}_{s\in S_{2r}}
\varepsilon_{q^{\pm 1}}(s) I^\pm_{k_{s(2)},k_{s(1)}} I^\pm_{k_{s(4)},k_{s(3)}}
 \cdots  I^\pm_{k_{s(2r)},k_{s(2r-1)}}.                           \label{Jpm}
\eeq
Here $1\le k_1<k_2<\ldots<k_{2r}\le n$, and the summation runs over
all the permutations $s$ of indices $k_1,k_2,\ldots,k_{2r}$ such that
\[
\begin{array}{c}
k_{s(2)}>k_{s(1)},\ \ \ \
k_{s(4)}>k_{s(3)},\ \ldots\ ,\ \  k_{s(2r)}>k_{s(2r-1)};\\
k_{s(2)}<k_{s(4)}<\ldots <k_{s(2r)}
\end{array}
\]
(the last chain of inequalities means that the sum includes only ordered
monomials). Symbol
$\varepsilon_{q^{\pm 1}}(s)\equiv (-q^{\pm 1})^{\ell(s)}$ stands for
a $q$-analogue of the Levi--Chivita antisymmetric tensor, $\ell(s)$ means
the length of permutation $s$.
(If $q\to 1$, both sets in (\ref{Jpm}) reduce
to the set of components of rank $2r$ antisymmetric tensor operator of the
Lie algebra ${\rm so}_n$.)

Using $q$-tensor operators given by (\ref{Jpm}) we obtain
the Casimir elements of $\Uqso$.

\ni {\bf Theorem 1.} {\it The  elements
\beq
C^{(2r)}_n=\sum_{1\le k_1<k_2<\ldots<k_{2r}\le n}
q^{k_1+k_2+\ldots+k_{2r}-r(n+1)}
J^+_{k_1,k_2,\ldots,k_{2r}} J^-_{k_1,k_2,\ldots,k_{2r}},      \label{GCE}
\eeq
where $r=1,2,\ldots, \bigl\lfloor \frac n2 \bigr\rfloor$
$(\lfloor x\rfloor$ means the integer part of $x)$, are Casimir elements
of $\Uqso$, i.e., they belong to the center of this algebra.

In fact, for even $n$, not only the product $($which constitutes
$C_n^{(n)}$\/$)$ of
elements $C^{(n)+}_n\equiv J^+_{1,2,\ldots,n}$ and
$C^{(n)-}_n\equiv J^-_{1,2,\ldots,n}$ belongs to the
center, but also each of them.
}

We conjecture that, in the case of $q$ being not a root of $1$, the set of
Casimir elements $C^{(2r)}_n$,
$r=1,2,\ldots, \bigl\lfloor \frac {n-1}2 \bigr\rfloor$, and
the Casimir element $C^{(n)+}_n$ (for even $n$) generates the center
of $\Uqso$, i.e., any element of the algebra $\Uqso$ which commutes
with all other elements can be presented as a polynomial of elements
from this set of Casimir elements.

Let us give explicitly some of Casimir elements. For $U'_q({\rm so}_3)$ and
$U'_q({\rm so}_4)$ we have
\[ %\hspace{-24mm}
C_3^{(2)}=q^{-1}I_{21}^2+I_{31}^+I_{31}^-+q I_{32}^2=
q I_{21}^2+I_{31}^-I_{31}^++q^{-1} I_{32}^2 ;
\]
\[ %\hspace{-23mm}
C_4^{(2)}=q^{-2}I_{21}^2+I_{32}^2+q^{2}I_{43}^2+
q^{-1}I_{31}^+I_{31}^-+q I_{42}^+I_{42}^-+I_{41}^+I_{41}^-,
\]
\[
C_4^{(4)+}=q^{-1}I_{21} I_{43} - I_{31}^+ I_{42}^+ + q I_{32}I_{41}^+=
q I_{21} I_{43} - I_{31}^- I_{42}^- + q^{-1} I_{32}I_{41}^-=C_4^{(4)-}.
\]
For $U'_q({\rm so}_5)$ the fourth order Casimir is
\[
C_5^{(4)}=q^{-2}J^+_{1,2,3,4} J^-_{1,2,3,4} + q^{-1} J^+_{1,2,3,5} J^-_{1,2,3,5}
+J^+_{1,2,4,5} J^-_{1,2,4,5} +
\]
\[
\hspace{0.5cm}+q J^+_{1,3,4,5} J^-_{1,3,4,5} +
q^2 J^+_{2,3,4,5} J^-_{2,3,4,5},
\]
where $J^+_{i,j,k,l}=q^{-1}I^+_{ji}I^+_{lk}-I^+_{ki}I^+_{lj}+q I^+_{kj}I^+_{li}$
and $J^-_{i,j,k,l}=q I^-_{ji}I^-_{lk}-I^-_{ki}I^-_{lj}+q^{-1} I^-_{kj}I^-_{li}$.
For $U'_q({\rm so}_6)$, we present only the highest order Casimir:
\begin{eqnarray}
C_6^{(6)+}
&=& q^{-3}I_{21}I_{43}I_{65}-q^{-2}I_{31}^+I_{42}^+I_{65} +
                                   q^{-1}I_{32}I_{41}^+I_{65}-\nonumber\\
&& -q^{-2}I_{21}I_{53}^+I_{64}^++q^{-1}I_{31}^+I_{52}^+I_{64}^+-
				I_{32}I_{51}^+I_{64}^++ \nonumber\\
&& +q^{-1}I_{21}I_{54}I_{63}^+-I_{41}^+I_{52}^+I_{63}^++
				q I_{42}^+I_{51}^+I_{63}^+- \nonumber\\
&& -I_{31}^+I_{54}I_{62}^++I_{41}^+I_{53}^+I_{62}^+-
				q^2 I_{43} I_{51}^+I_{62}^++ \nonumber\\
&& +q I_{32}I_{54}I_{61}^+-q^2 I_{42}^+I_{53}^+I_{61}^++
				q^3 I_{43} I_{52}^+I_{61}^+. \nonumber
\end{eqnarray}
Finally, let us give explicitly the quadratic Casimir of $\Uqso$,
\[
C_n^{(2)}=\sum_{1\le i < j \le n} q^{i+j-n-1} I_{ji}^+ I_{ji}^-.
\]
This formula coincides with that
given in \cite{NUW}, and is a particular case of (\ref{GCE}).

\medskip
\ni {\bf 4. Irreducible representations of $\Uqso$.}
\medskip

Let us give a brief description of irreducible representations (irreps) of
$\Uqso$. More detailed description of these irreps can be found in
\cite{GK,GI}.

As in the case of Lie algebra ${\rm so}_n$,
finite-dimensional irreps $T$
of the algebra $\Uqso$ are characterized by the set ${\bf m}_n\equiv$
$(m_{1,n},$ $m_{2,n},\ldots,$ $m_{\lfloor\frac n2 \rfloor,n})$
(here $\lfloor x\rfloor$ means the integer part of $x$) of numbers, which
are either all integers or all half-integers, and satisfy the well-known
dominance conditions
\[
\begin{array}{ll}
m_{1,n}\ge m_{2,n}\ge\cdots\ge m_{\frac n2 -1,n}\ge
|m_{\frac n2, n}|\ \ \ \
& \mbox{if\  $n$\  is even,} \\
m_{1,n}\ge m_{2,n}\ge\cdots\ge m_{\frac {n-1}2,n}\ge 0 &
\mbox{if\  $n$\  is odd.}
\end{array}
\]
To give the representations in Gel'fand--Tsetlin basis we denote, as in
the case of Lie algebra ${\rm so}_n$, the basis vectors $|\alpha\rangle$
of representation
spaces by Gel'fand--Tsetlin patterns $\alpha$.
The representation operators
$T_{{\bf m}_n}(I_{2p+1,2p})$ and $T_{{\bf m}_n}(I_{2p,2p-1})$ act on
$|\alpha\rangle$ by the formulas
\[ %\hspace{-20mm}
T_{{\bf m}_n}(I_{2p+1,2p})|\alpha\rangle=
\sum_{r=1}^p \Bigl(A^r_{2p}(\alpha)|m^{+r}_{2p}\rangle-
A^r_{2p}(m^{-r}_{2p})|m^{-r}_{2p}\rangle\Bigr),
\]

\vspace{-3mm}
\[
T_{{\bf m}_n}(I_{2p,2p-1})|\alpha\rangle=
\sum_{r=1}^{p-1} \Bigl(B^r_{2p-1}(\alpha)|m^{+r}_{2p-1}\rangle-
B^r_{2p-1}(m^{-r}_{2p-1})|m^{-r}_{2p-1}\rangle\Bigr)+
{\rm i}C_{2p-1}|\alpha\rangle.
\]
Here the matrix elements $A^r_{2p}$,$B^r_{2p-1}$,
$C_{2p-1}$ are obtained from the classical (non-deformed)
ones by replacing each
factor $(x)$ with its respective $q$-number
$[x]\equiv(q^x-q^{-x})/(q-q^{-1})$; besides,
the coefficient $\frac12$ in the `classical' $A^r_{2p}$ is replaced with
the $l_{r,2p}$-dependent expression
$\bigl(([l_{r,2p}][l_{r,2p}+1])/ ([2l_{r,2p}][2l_{r,2p}+2]) \bigr)^
{\hspace{-0.2mm}1/2}$,
where $l_{r,2p}=m_{r,2p}+p-r$.

\medskip
\ni {\bf 5. Casimir operators and their eigenvalues.}
\medskip

The Casimir operators
(the operators which correspond to the Casimir elements), within
irreducible finite-dimensional representations of $\Uqso$ take
diagonal form.
To give them explicitly,
we employ the so-called {\it generalized factorial
elementary symmetric polynomials} (see \cite{MN}). Fix an arbitrary
sequence of complex numbers ${\bf a}=(a_1,a_2,\ldots)$. Then, for each
$r=0,1,2,\ldots,N$, introduce the polynomials
of $N$ variables $z_1,$ $z_2,\ldots,$ $z_N$ as follows:
\beq
e_r(z_1,z_2,\ldots,z_N|{\bf a})=\sum_{1\le p_1<p_2<\cdots<p_r\le N}
(z_{p_1}-a_{p_1}) (z_{p_2}-a_{p_2-1})\ldots (z_{p_r}-a_{p_r-r+1}).
\eeq

The Casimir operators in the
irreducible finite-dimensional representations characterized by the set
$(m_{1,n},$ $m_{2,n},\ldots,$ $m_{N,n})$,
$N=\bigl\lfloor\frac n2\bigr\rfloor$, by the Schur Lemma,
are presentable as (here {\bf 1} denotes the unit operator):
\[
T_{{\bf m}_n}(C^{(2r)}_n)=\chi^{(2r)}_{{\bf m}_n}{\bf 1}.
\]

\ni {\bf Theorem 2.} {\it The eigenvalue of the operator
$T_{{\bf m}_n}(C^{(2r)}_n)$ is
\[
\chi^{(2r)}_{{\bf m}_n}=(-1)^r
e_r([l_{1,n}]^2,[l_{2,n}]^2,\ldots,[l_{N,n}]^2|{\bf a})
\]
where ${\bf a}=([\epsilon]^2,[\epsilon+1]^2,[\epsilon+2]^2,\ldots)$,
$l_{k,n}=m_{k,n}+N-k+\epsilon$. Here $\epsilon=0$ for $n=2N$ and
$\epsilon=\frac12$ for $n=2N+1$.

In the case of even $n$, i.e., $n=2N$,
\[
T_{{\bf m}_n}(C^{(n)+}_n)=T_{{\bf m}_n}(C^{(n)-}_n)=\bigl(\sqrt{-1}\bigr)^N
[l_{1,n}][l_{2,n}]\ldots[l_{N,n}] {\bf 1}.
\]
}

The eigenvalues of Casimir operators are important
for physical applications.
Let us quote some of Casimir operators together with their eigenvalues.
For $U'_q({\rm so}_3)$,
\[
T_{(m_{13})}(C^{(2)}_3)=-[m_{13}][m_{13}+1] {\bf 1}.
\]
For $U'_q({\rm so}_4)$ we have
\[ %\hspace{-12mm}
T_{(m_{14},m_{24})}(C^{(2)}_4)=-([m_{14}+1]^2+[m_{24}]^2-1) {\bf 1},
\]
\[
T_{(m_{14},m_{24})}(C^{(4)+}_4)=T_{(m_{14},m_{24})}(C^{(4)-}_4)=
-[m_{14}+1][m_{24}]{\bf 1}.
\]
Finally, for $U'_q({\rm so}_5)$ the Casimir operators are
\[
T_{(m_{15},m_{25})}(C^{(2)}_5)=-([m_{15}+3/2]^2+[m_{25}+1/2]^2-[1/2]^2-
                      [3/2]^2) {\bf 1},
\]
\[
T_{(m_{15},m_{25})}(C^{(4)}_5)=([m_{15}+3/2]^2-[1/2]^2)([m_{25}+1/2]^2-[1/2]^2)
                                                {\bf 1}.
\]

\medskip
\ni {\bf 6. Concluding remarks.}
\medskip

In this note, for the nonstandard $q$-algebras $\Uqso$
we have presented explicit formulae for all the Casimir
operators corresponding to basis set of Casimirs of ${\rm so}_n$.
Their eigenvalues in irreducible finite-dimensional representations are
also given.
We believe that the described Casimir elements
generate the whole center of the algebra $\Uqso$
(of course, for $q$ being not a root of unity).

As mentioned, the algebras $\Uqso$ for $n>4$ are of importance
in the {\it construction of algebra of observables for} 2+1
{\it quantum gravity}
(with 2D space of genus $g>1$) serving as certain intermediate algebras.
For that reason, the results concerning
Casimir operators and their eigenvalues will be useful in the process
of construction of
the desired algebra of independent quantum observables for the case of
higher genus surfaces, in the important and interesting case of
anti-De Sitter gravity (corresponding to negative
cosmological constant).

%%%%%%%%%%%%%%%%%%%%%%%%
\smallskip
\ni {\bf Acknowledgements.}
The authors express their gratitude to J.~Nelson and
A.~Klimyk for interesting discussions.
This research was supported in part by the CRDF Award No. UP1-309, and by
the DFFD Grant 1.4/206.

\end{document}